\documentclass[a4paper,twoside,11pt,righttag]{amsart}
\usepackage{amsxtra}
\usepackage{amsfonts}
\usepackage{amssymb, verbatim}
\usepackage{graphics,epsfig}

\theoremstyle{plain}
\newtheorem{prop}{Proposition}[section]
\newtheorem{teo}[prop]{Theorem}
\newtheorem{coro}[prop]{Corollary}

\theoremstyle{definition}

\newtheorem{ejem}[prop]{Example}
\newtheorem{rem}[prop]{Remark}

\theoremstyle{remark}

\numberwithin{equation}{section}

\newcommand{\Z}{\mathbb Z}

\newcommand{\R}{\mathbb R}

\newcommand{\g}{\gamma}
\newcommand{\ld}{\lambda}
\newcommand{\Ld}{\Lambda}
\newcommand{\G}{\Gamma}
\newcommand{\D}{\Delta}

\newcommand{\Lf}{\D_{f}}

\newcommand{\vcp}{\Gamma \backslash \R^n}
\newcommand{\I}{\text{\sl Id}}

\newcommand{\man}{M_\Gamma}
\newcommand{\tor}{T_\Ld}

\newcommand{\on}{\text{O}(n)}

\newcommand{\fa}{\mathcal{F}}

\newcommand{\bjh}{B_{j,h}}

\title[$\Z_2^k$-manifolds are isospectral on forms.]
{$\Z_2^k$-MANIFOLDS ARE ISOSPECTRAL ON FORMS.}
\author[R. J. Miatello]{R. J. Miatello}
\address{FaMAF--CIEM \\$\!$Universidad Nacional de C\'ordoba\\5000 C\'ordoba, Argentina.}
\email{miatello@mate.uncor.edu}
\author[R. A. Podest\'a]{R. A. Podest\'a}
\email{podesta@mate.uncor.edu}
\author[J. P. Rossetti]{J. P. Rossetti.}
\email{rossetti@mate.uncor.edu}
\thanks{Supported by Conicet, Secyt-UNC}

\begin{document}
\bibliographystyle{plain}

\begin{abstract}
We obtain a simple formula for the multiplicity of eigenvalues of
the Hodge-Laplace operator, $\Lf$, acting  on sections of the full
exterior bundle $\Lambda(TM)=\bigoplus_{p=0}^n \Lambda^p(TM)$ over
an arbitrary compact flat Riemannian $n$-manifold $M$ with
holonomy group $\Z_2^k$, with $1 \le k \le n-1$. This formula
implies that any two compact flat manifolds with holonomy group
$\Z_2^k$ having isospectral lattices of translations are
isospectral on forms, that is, with respect to $\Lf$. As a
consequence, we construct a large family of pairwise
$\Lf$-isospectral and nonhomeomorphic $n$-manifolds of cardinality
greater than $2^{\frac{(n-1)(n-2)}2}$.
\end{abstract}

\maketitle

\section*{Introduction}     \label{s.intro}
In [{\bf MR2,3,4}] the spectrum of the Hodge-Laplacian on
$p$-forms on compact flat manifolds was studied, comparing
$p$-isospectrality with other types of isospectrality. In
particular, pairs of manifolds that are isospectral on $p$-forms
for a fixed value of $p>0$ were constructed,  having different
lengths of closed geodesics or different first eigenvalue of the
Laplacian on functions. Most of the examples given belong to the
class of {\em $\Z_2^k$-manifolds}, that is, flat Riemannian
manifolds with holonomy group $\Z_2^k$. By the
Cartan-Ambrose-Singer theorem, such manifolds are necessarily
flat, hence of the form $M_\G=\G\backslash \R^n$, $\G$ a
Bieberbach group with translation lattice $\Lambda$ and with
holonomy group $F:=\Lambda \backslash \G\simeq \Z_2^k$.

The goal of this paper is to show that if we consider the {\em
full}  exterior bundle over a  $\Z_2^k$-manifold $\man$,  there is
a high degree of regularity in the spectrum of the Hodge
Laplacian, $\Lf$, acting on sections of this bundle. Two manifolds
having the same spectrum with respect to $\Lf$ will be called {\em
isospectral on forms}. We shall see that the spectrum of a flat
manifold $\man$ is completely determined by the spectrum of the
covering torus $T_\Ld$, and furthermore any two $\Z_2^k$-manifolds
$\man$, $M_{\G'}$, with covering torus $T_\Ld$ are isospectral on
forms. They are also isospectral on even (resp.\@ odd) forms, that
is, with respect to the operator $\Delta_f$ restricted to even
(resp.\@ odd) forms. This allows to obtain very large families of
$\Lf$-isospectral $n$-manifolds, pairwise nonhomeomorphic to each
other. In particular we will describe a family of flat manifolds,
the so called generalized Hantzsche-Wendt manifolds (see
\cite{RS}), having holonomy group $\Z_2^{n-1}$, whose cardinality
is greater than $2^{\frac{(n-1)(n-2)}2}$. The proof of the main
result uses the multiplicity formulae in \cite{MRp} together with
some symmetry properties of the Krawtchouk polynomials. We point
out that the above isospectrality result is valid only for
holonomy groups $F\simeq \Z_2^k$. Indeed, we shall see  that it
fails to hold for flat manifolds with holonomy group $F\simeq
\Z_4$ and $F\simeq \Z_4\times \Z_2$ (Example 3.5).

\section{Preliminaries}
We first recall from \cite{Ch} or  \cite{Wo} some standard facts
on compact flat manifolds. A Bieberbach  group is a discrete,
cocompact torsion-free subgroup  $\G$ of the isometry group
$I(\R^n)$ of $\R^n$. Such  $\G$ acts properly discontinuously on
$\R^n$, hence $M_\G = \vcp$ is   a compact flat Riemannian
manifold with fundamental group $\G$. Furthermore, any such
manifold arises in this way. Since $I(\R^n)\simeq \on \ltimes
\R^n$, any element $\g \in I(\R^n)$ decomposes uniquely as $\g = B
L_b$, with $B \in \on$ and $b\in \R^n.$ The translations in $\G$
form a normal maximal abelian subgroup  of finite index $L_\Ld$,
$\Ld$ a lattice in $\R^n$ which is $B$-stable for every $BL_b \in
\G$. The restriction to $\G$ of the canonical projection
$r:I(\R^n) \rightarrow \on$, given by $BL_b \mapsto B$, is a
homomorphism with kernel $\Ld$ and $r(\G)$ is a finite subgroup of
$\on$ isomorphic to $F:= \Ld \backslash \G$, the linear holonomy
group of the Riemannian manifold $M_\G$.

We recall from \cite{MRp} the multiplicity formula for the
eigenvalues of the Hodge Laplace operator $-\D_p$ acting on smooth
$p$-forms of a compact flat manifold $\man$. For any $\mu \ge 0$,
let
\begin{equation}
\Ld^*_\mu=\{\ld \in \Ld^* : \| \ld\|^2=\mu\}
\end{equation}
 where $\Ld^\ast$ is the dual lattice of $\Ld$. In \cite{MRp}, Theorem 3.1, it is
shown that the multiplicity of the eigenvalue $4\pi^2 \mu$ of $-\Delta_p$
is given by
\begin{equation}                                \label{dpmu}
 d_{p,\mu}(\G)= \tfrac{1}{|F|} \sum_{\g=BL_b
 \in \Ld\backslash \G}\text{tr}_p(B) \; e_{\mu,\g}(\G)
\end{equation}
where $e_{\mu,\g} = \sum_{v\in {\Ld^*_\mu}:Bv=v} e^{-2\pi i
 v\cdot b}$ and $\text{tr}_p$ is the trace of the $p$-exterior
 representation $\tau_p:\on \rightarrow \text{GL}(\Lambda^p(\R^n))$.

A Bieberbach group $\G$ is said to be of {\em diagonal type} (see
\cite{MRd}, Definition 1.3) if there exists an orthonormal
$\Z$-basis $\{\ld_1,\dots,\ld_n\}$ of the lattice $\Ld$ such that
for any element $BL_b\in\G$, $B\ld_i=\pm \ld_i$ for $1\le i\le n$.
These Bieberbach groups have holonomy group $\Z_2^k$ for some
$1\le k\le n-1$. If $\G$ is of diagonal type, after conjugation of
$\G$ by an isometry,  it may be assumed that $\Ld$ is the
canonical (or cubic)  lattice and, furthermore, that $b$ lies in
$\frac 1{2} \Ld$ for any $\g=BL_b \in \G$ (see \cite{MRd}, Lemma
1.4).

For Bieberbach groups of diagonal type, the traces
$\text{tr}_p(B)$ in (\ref{dpmu}) are given by integral values of
the {\em Krawtchouk polynomials} of degree $p$
\begin{equation} \label{kraw}
K_p^n(x):= \sum_{t=0}^p (-1)^t \binom xt \binom
    {n-x}{p-t}
\end{equation}
(see \cite{MRp}, Remark~3.6 and \cite{MRd}; also, see \cite{KL} for more
information on Krawtchouk polynomials).
 Indeed, we have
 \begin{equation} \label{tracep}
    \text{tr}_p(B)=K_p^n(n-n_B),
\quad \text{where } n_B:= \dim \, (\R^n)^B = \dim \ker(B-\I).
\end{equation}

The first Krawtchouk polynomials are $K^n_0(x)=1$, $K^n_1(x)=-2x
+n$, $K^n_2(x)=2x^2-2nx+\tbinom n2$,  $K^n_3(x)=-\frac43 x^3+
2nx^2 -(n^2-n+\frac23)x + \tbinom n3$. For later use we also give,
in the following tables, the integral values of $K^n_p(x)$ for
$0\le p,x \le n$,  $n=3,4$.

\renewcommand{\arraystretch}{0.2}
\begin{center}
\begin{equation} \label{krawtables}
\begin{tabular}{|c|cccc|}  \hline & & & &  \\
$x$ & 0 & 1 & 2 & 3  \\ & & & & \\  \hline & & & & \\
$K_0^3(x)$ & 1 & 1 & 1 & 1 \\ & & & & \\
$K_1^3(x)$ & 3 & 1 & -1& -3  \\ & & & & \\
$K_2^3(x)$ & 3 & -1& -1 & 3  \\ & & & & \\
$K_3^3(x)$ & 1 & -1 & 1 & -1  \\  & & & & \\
 \hline
\end{tabular} \qquad
\begin{tabular}{|c|ccccc|}  \hline & & & & &  \\
$x$ & 0 & 1 & 2 & 3 & 4  \\ & & & & & \\  \hline & & & & & \\
$K_0^4(x)$ & 1 & 1 & 1 & 1 & 1 \\ & & & & & \\
$K_1^4(x)$ & 4 & 2 & 0 & -2 & -4 \\ & & & & & \\
$K_2^4(x)$ & 6 & 0 & -2 & 0 & 6 \\ & & & & & \\
$K_3^4(x)$ & 4 & -2 & 0 & 2 & -4 \\ & & & & & \\
$K_4^4(x)$ & 1 & -1 & 1 & -1 & 1 \\ & & & & & \\ \hline
\end{tabular}
\end{equation}
\end{center}

\section{The spectrum on forms of $\Z_2^k$-manifolds.}
Let $\bigoplus_{p=0}^n \bigwedge^p(T(\man))$ be the full exterior
bundle of the compact flat manifold $\man$ and let $\Delta_p$ be
the Hodge Laplacian acting on $p$-forms. We shall denote by
\begin{equation} \label{formlaplacian}
\D_{f} := \sum_{p=0}^n \D_p, \qquad \D_{e} := \sum_{p \text{
even}} \D_p, \qquad \D_{o} := \sum_{p \text{ odd}} \D_p,
\end{equation}
 the {\em Laplacian on forms},  {\em on even forms}   and {\em  on odd forms} of $M_\G$,
 respectively.

The multiplicity of the eigenvalue $4\pi^2\mu$ for $\Lf$ is given
by
\begin{equation}\label{dfmu}
 d_{f,\mu}(\G)=\sum_{p=0}^n d_{p,\mu}(\G)
\end{equation}
and similarly   $d_{e,\mu}(\G)=\sum_{p \text{ even}}
d_{p,\mu}(\G)$ and $d_{o,\mu}(\G) =\sum_{p \text{ odd}}
d_{p,\mu}(\G)$, for $\D_{e}$ and  $\D_{o}$ respectively. Thus,
$\Lf=\D_{e}+\D_{o}$ and
$d_{f,\mu}(\G)=d_{e,\mu}(\G)+d_{o,\mu}(\G)$.

Clearly, $p$-isospectrality for all $p$ implies $\Lf$-isospectrality (as
well as $\D_{e}$ and $\D_{o}$-isospectrality), but
 we shall see that the converse is far from being true.

\begin{teo} \label{main}
If $\G$ is a Bieberbach group with translation lattice $\Ld$ and
holonomy group $\Z_2^k$, then for any $\mu \ge 0$  the
multiplicities of the eigenvalue $4\pi^2\mu$ for $\Lf$,  $\D_{e}$
and $\D_{o}$ are given respectively by
$$d_{f,\mu}(\G) = 2^{n-k}|\Ld^\ast_\mu|,  \qquad
d_{e,\mu}(\G)= d_{o,\mu}(\G) = 2^{n-k-1}|\Ld^\ast_\mu|.$$ Thus, if
$M_\G, M_{\G'}$ are $\Z_2^k$-manifolds with translation lattices
$\Ld, \Ld'$, then $M_\G$ and $M_{\G'}$ are isospectral on forms
(resp.\@ on even or odd forms) if and only if $\Ld$ and $\Ld'$ are
isospectral. In particular, for fixed $\Ld$ and $k$, {\em all}
$\Z_2^k$-manifolds having covering torus $\tor$ are $\Lf$,
$\D_{e}$ and $\D_{o}$-isospectral.
\end{teo}

\begin{proof}
Let $\man$ be a $\Z_2^k$-manifold. Then  $M_\G = \vcp$ with
$\G=\langle \g_1,\dots,\g_k, \Ld \rangle$, where  $\Ld$ is a
lattice and $\g_i=B_iL_{b_i}$, $B_i \in \on$, $b_i \in \R^n$,
$B_i\Ld=\Ld$, $B_i^2=\I$,  $B_i B_j=B_j B_i$ for each $1\leq
i,j\leq k$.

We know that if $B_i$ is diagonal then $\text{tr}_p \,(B_i) =
K_p^n(n-n_{B_i})$ (see \cite{MRp}, Remark 3.6). This is also valid
for non-diagonal matrices $B$ of order 2. Indeed, $B$ has only
eigenvalues of the form $\pm 1$, hence $B$ is conjugate in
$\text{GL}_n(\R)$ to the diagonal matrix $D_{B}:=\left[
\begin{smallmatrix} -\text{I}_{n-n_B} & \\ & \text{I}_{n_B}
\end{smallmatrix} \right]$ where $\text{I}_m$ is the identity
matrix in $\R^m$. Thus $\text{tr}_p(B)=\text{tr}_p(D_B)=K_p^n(n-n_B)$.

Hence,
 by (\ref{dpmu}), (\ref{tracep}) and the fact that $K_p^n(0) =\tbinom
 np$, we have
\begin{eqnarray*}
d_{p,\mu}(\G)
  =  2^{-k} \Big(\tbinom np |\Ld^\ast_\mu| + \sum_{\g \in \Ld\backslash \G,\,\, \g \ne \I} K_p^n(n-n_B)
    \; e_{\mu,\g}(\G)\Big)
 \end{eqnarray*}
and, adding over $p$, we obtain
\begin{equation*}
d_{f,\mu}(\G)
 = 2^{n-k}|\Ld^\ast_\mu| +
 2^{-k} \sum_{\g \in \Ld\backslash \G, \,\, \g \ne \I} \Big( \sum_{p=0}^n K_p^n(n-n_B)
 \Big) \; e_{\mu,\g}(\G).
 \end{equation*}

Now, we show that $\sum_{p=0}^n K_p^n(j)=0$ for fixed $j\not=0$.
In fact,
\begin{eqnarray*} \label{krawsum}
\sum_{p=0}^n K_p^n(j) & = & \sum_{p=0}^n \sum_{t=0}^j (-1)^t
\tbinom jt \tbinom{n-j}{p-t}\\ & =& \sum_{t=0}^j (-1)^t \tbinom jt
\sum_{p-t=0}^{n-j} \tbinom{n-j}{p-t}  \\&=&   2^{n-j} \sum_{t=0}^j
(-1)^t \tbinom jt =0
\end{eqnarray*}
Thus, since $n-n_B=0$ if and only if $B=\I$, we obtain that
$d_{f,\mu}(\G)= 2^{n-k} |\Ld^\ast_\mu|$, as claimed. The proofs
for $d_{e,\mu}(\G)$ and $d_{o,\mu}(\G)$ are the same, except that
we add  over even and odd values of $p$, respectively.
\end{proof}

\begin{rem}
(i) We note that for the $n$-torus $\tor$, for each $0\le p\le n$,
we have $d_{p,\mu}(\Ld)=\tbinom np |\Ld_\mu^\ast|$, hence
$0$-isospectrality is equivalent to $p$-isospec\-tra\-lity for any
$p>0$, and this in turn is equivalent to $\Lf$-isospectrality.
However, there are many examples of pairs of compact flat
manifolds that are $p$-isospectral for some $p>0$ but are not
isospectral on functions and also pairs of manifolds that are
isospectral on functions  but are not $p$-isospectral for any
$0<p<n$ (see \cite{MRp, MRd}).

(ii) We shall see that Theorem \ref{main} does not hold for general holonomy
groups. For instance, Example 3.5 will show it  fails
to hold when $F$ is isomorphic  to $\Z_4$ or $\Z_4 \times \Z_2$.
\end{rem}

\section{Examples and Counterexamples}

\begin{ejem} We now consider a
family of $\Z_2$-manifolds, of cardinality quadratic in $n$, which
are pairwise not isospectral on functions, but which are
isospectral on forms, according to Theorem \ref{main}.

Put $J:=\left[\begin{smallmatrix} 0 & 1 \\
1 & 0
\end{smallmatrix} \right]$. For each $0\leq j,h < n$, define
\begin{equation}
 B_{j,h}:=\text{diag}(\underbrace{J,\dots,J}_j,
 \underbrace{-1,\dots,-1}_h,\underbrace{1,\dots,1}_l)
\end{equation}
 where $n=2j+h+l$, $j+h\not=0$ and $l\geq1$.  Then $\bjh \in \on$,
 $\bjh^2=\I$.
Let $\Ld=\Z e_1\oplus \cdots \oplus \Z e_n$ be the canonical lattice of
 $\R^n$ and for $j,h$ as before define the groups
\begin{equation}
 \G_{j,h}:=\langle \bjh L_{\frac{e_n}{2}}, \Ld \rangle.
\end{equation}
 We have that $\Ld$ is stable by $B_{j,h}$ and
 $(\bjh +\I)\frac{e_n}2 = e_n \in \Ld \smallsetminus (\bjh +\I)\Ld$.
 It is easy to verify that $\G_{j,h}$ is torsion-free, hence a  Bieberbach group.
   In this way, if we set $M_{j,h}:= \G_{j,h} \backslash \R^n$, we have
     a family
\begin{equation}                                        \label{family}
 \mathcal{F}:=\{ M_{j,h}=\G_{j,h} \backslash \R^n \,:\, 0 \le j \le [\tfrac{n-1}2],
 0 \le h < n-2j, \,j+h\not=0 \}
\end{equation}
of compact flat manifolds  with holonomy group $F\simeq \Z_2$.

Furthermore, the family $\mathcal{F}$ gives a  system of
representatives for the diffeomorphism classes of $\Z_2$-manifolds
of dimension $n$ (see \cite{MP} for a proof). Also
\begin{equation}                                        \label{H1}
 H_1(M_{j,h},\Z) \simeq \Z^{j+l} \oplus \Z_2^h,
\end{equation}
and if $1\le p\le n$, then the Betti numbers are
\begin{equation}                                        \label{bettip}
 \beta_p(M_{j,h})=  \sum _{i=0}^{[\frac p2]}\binom{j+h}{2i}
 \binom{j+l}{p-2i}.
\end{equation}
Hence, if $\beta_1(M_{j,h})=\beta_1(M_{j',h'})$, then
$\beta_p(M_{j,h})=\beta_p(M_{j',h'})$ for any $p \ge 1$.

 Moreover,
\begin{equation}                                        \label{e.card}
 \#\fa  =  \big(n-[\tfrac{n-1}2]\big)\big([\tfrac{n-1}2]+1\big)-1
 =\left\{ \begin{array}{lr}\frac{n^2+2n-4}{4} & \qquad \text{$n$ even} \\
 \frac{n^2+2n-3}{4} & \qquad \text{$n$ odd.} \end{array} \right.
\end{equation}

Now, since $\bjh$ and $B_{0,j+h}$ are conjugate in
$\text{GL}_n(\R)$ we have that $\text{tr}_p(\bjh)$ $=
\text{tr}_p(B_{0,j+h}) = K_p^n(j+h)$. Hence, by formula
(\ref{dpmu}), the expression for the multiplicity of the
eigenvalue $4\pi^2\mu$ of $-\Delta_p$ equals
\begin{equation}    \label{e.z2mult}
d_{p,\mu}(\G_{j,h})=\tfrac12 \big( \tbinom{n}{p} |\Ld_\mu| + K_p^n(j+h)
\,e_{\mu,\g}(\G_{j,h}) \big).
\end{equation}
where $e_{\mu,\g}(\G) = \sum_{v\in \Ld^B_\mu} e^{-2\pi i v\cdot
b}$ and $\Ld^B$ denotes the  elements in $\Ld$ fixed by~$B$.

We claim that the manifolds in $\fa$ are pairwise not isospectral
on functions. To see this, it will suffice to compare the
multiplicities of the two smallest nonzero eigenvalues, namely
$\mu=1$ and $\mu=\sqrt 2$.

Take $\mu=1$. Then $\Ld_1=\{\pm e_1,\dots,\pm e_n\}$ and
$\Ld_1^{\bjh} = \{\pm e_{2j+h+1},\dots,\pm e_n\}$ thus
$|\Ld_1|=2n$ and $|\Ld_1^{\bjh}|=2(n-(2j+h))=2l$. Now, one checks
that $e_{1,\g}(\G_{j,h})=2(l-1)+2(-1)=2(l-2)$ and hence we get
from (\ref{e.z2mult})
\begin{equation} \label{dp1}
d_{p,1}(\G_{j,h})= \tbinom{n}{p} n + K^n_p(j+h)(l-2).
\end{equation}

Now consider $\mu=\sqrt 2$. Then
 $\Ld_{\sqrt 2}=\{\pm (e_i \pm e_j) : 1\le i < j \le n\}$ and
 $\Ld_{\sqrt 2}^{B_{j,h}}=\{\pm (e_{2i-1} + e_{2i}) : 1\le i \le j\}\cup
 \{ \pm(e_i\pm e_j) : n-l+1\le i < j \le n\}$. Hence
 $|\Ld_{\sqrt 2}|=4\tbinom n2$ and $|\Ld_{\sqrt 2}^{\bjh}|=2j+4 \tbinom l2$.
One checks that $e_{\sqrt 2,\g}(\G_{j,h})=2j + 4\left(\binom l2
-2(l-1)\right)=2j+2(l-1)(l-4)$. In this way we obtain
\begin{equation} \label{dp2}
d_{p,\sqrt 2}(\G_{j,h})= 2\tbinom{n}{p} \tbinom n2 +
K^n_p(j+h)\big(j+(l-1)(l-4)\big).
\end{equation}

In particular for $p=0$, since $K^n_0(j)=1$ for any $j$, we have
\begin{eqnarray}\label{d1}
d_{0,1}(\G_{j,h})& = & n + l - 2 \\
\label{d2} d_{0,\sqrt 2}(\G_{j,h})& =& n(n-1) + j + (l-1)(l-4).
\end{eqnarray}

This allows to distinguish the spectra on functions of the
$\Z_2$-manifolds considered. Indeed, if $M_{j,h}, M_{j',h'}$ are
isospectral then $l=l'$ by (\ref{d1}), thus $2j+h=2j'+h'$. By
(\ref{d2}), $j=j'$ and hence $h=h'$. This shows that all manifolds
in $\mathcal{F}$ are pairwise not isospectral to each other.

To illustrate the compensations occurring in the sums  in
(\ref{dfmu}) we compute the individual multiplicities
$d_{p,\mu}(\G)$ corresponding to $\mu=1,\sqrt 2$, for manifolds in
$\mathcal F$ in dimensions 3 and 4.

In dimension 3, there are only three $\Z_2$-manifolds up to
diffeomorphism (see \cite{Wo}): $M_{1,0}, M_{0,2}$ and
$M_{0,1}$, with holonomy groups generated respectively
by the matrices $[\begin{smallmatrix} J &  \\
& 1 \end{smallmatrix}]$,  $[\begin{smallmatrix} -I & \\ & 1
\end{smallmatrix}]$ and $[\begin{smallmatrix} -1 &  \\ & I
\end{smallmatrix}]$ where $I$ is the $2\times 2$ identity
matrix. These manifolds are called  {\em dicosm} ($c2$), {\em
first amphicosm} ($+a1$)  and {\em second amphicosm} ($-a1$)
respectively, in \cite{CR}.

Using formulae (\ref{dp1}) and (\ref{dp2}) and the tables in
(\ref{krawtables}) for the integral values of Krawtchouk
polynomials we compute the following values of $d_{p,1}$ and
$d_{p, \sqrt 2}$, for $0 \le p \le 3$:

\renewcommand{\arraystretch}{0.3}
\begin{center}
\begin{tabular}{|c|cccc|c|}  \hline & & & & & \\
$\mu=1$ & $d_0$ & $d_1$ & $d_2$ & $d_3$ & $d_{f}$ \\
& & & & & \\  \hline & & & & & \\
$M_{1,0}$ & 2 & 8 & 10 & 4 & 24 \\ & & & & & \\
$M_{0,2}$ & 2 & 10& 10 & 2 & 24 \\ & & & & & \\
$M_{0,1}$ & 3 & 9 & 9  & 3 & 24 \\ & & & & & \\ \hline
\end{tabular}\qquad
\renewcommand{\arraystretch}{0.27}
\begin{tabular}{|c|cccc|c|}  \hline & & & & & \\
$\mu=\sqrt2$ & $d_0$ & $d_1$ & $d_2$ & $d_3$ & $d_{f}$ \\
& & & & & \\  \hline & & & & & \\
$M_{1,0}$ & 7 & 19 & 17 & 5 & 48 \\ & & & & & \\
$M_{0,2}$ & 6 & 18& 18 & 6 & 48 \\ & & & & & \\
$M_{0,1}$ & 4 & 16 & 20  & 8 & 48 \\ & & & & & \\ \hline
\end{tabular}
\end{center}

\medskip
In dimension 4 there are five nondiffeomorphic $\Z_2$-manifolds,
$M_{1,1}, M_{1,0}$, $M_{0,3}$, $M_{0,2}$ and $M_{0,1}$, with
holonomy group generated, respectively, by the matrices $\left[\begin{smallmatrix} J & & \\
& -1 &
\\ & & 1
\end{smallmatrix}\right]$, $\left[\begin{smallmatrix} J & & \\ & 1 & \\ & & 1
\end{smallmatrix}\right]$, $\left[\begin{smallmatrix} -I & & \\ & -1 & \\ & & 1
\end{smallmatrix}\right]$, $\left[\begin{smallmatrix} -I & & \\ & 1 & \\ & & 1
\end{smallmatrix}\right]$ and $\left[\begin{smallmatrix} -1 & & \\ & 1 & \\ & & I
\end{smallmatrix}\right]$. Proceeding as before we get the tables:
\renewcommand{\arraystretch}{0.3}
\begin{center}
\begin{tabular}{|c|ccccc|c|}  \hline & & & & & & \\
$\mu=1$ & $d_0$ & $d_1$ & $d_2$ & $d_3$ & $d_4$&  $d_{f}$ \\
& & & & & &\\  \hline & & & & & & \\
$M_{1,1}$ & 3 & 16& 26 & 16& 3& 64 \\ & & & & & & \\
$M_{1,0}$ & 4 & 16& 24 & 16& 4& 64 \\ & & & & & & \\
$M_{0,3}$ & 3 & 18& 24 & 14& 5& 64 \\ & & & & & & \\
$M_{0,2}$ & 4 & 16& 24 & 16& 4& 64 \\ & & & & & & \\
$M_{0,1}$ & 5 & 18& 24 & 14& 3& 64 \\ & & & & & & \\ \hline
\end{tabular}\quad
\renewcommand{\arraystretch}{0.28}
\begin{tabular}{|c|ccccc|c|}  \hline & & & & & & \\
$\mu=\sqrt2 $ & $d_0$ & $d_1$ & $d_2$ & $d_3$ & $d_4$&  $d_{f}$ \\
& & & & & &\\  \hline & & & & & & \\
$M_{1,1}$ & 13& 48& 70 & 48&13& 192\\ & & & & & & \\
$M_{1,0}$ & 11& 46& 72 & 50&13& 192\\ & & & & & & \\
$M_{0,3}$ & 12& 48& 72 & 48&12& 192\\ & & & & & & \\
$M_{0,2}$ & 10& 48& 76 & 48&10& 192\\ & & & & & & \\
$M_{0,1}$ & 10& 44& 72 & 52&14& 192\\ & & & & & & \\ \hline
\end{tabular}
\end{center}
\end{ejem}

\medskip
\begin{ejem}\label{ej3.2}
Here we consider the $\Z_2^2$-manifolds of dimension 3 having the
cubic lattice as lattice of translations (see \cite{Wo}, Section
3.5). There are three such manifolds, up to isometry. We shall see
that they are not $p$-isospectral for any $0\le p\le n$, showing
for small eigenvalues how the compensations take place so that the
sums of multiplicities for all $p$ become the same.

We consider the Hantzsche-Wendt manifold, $M_1=\G_1\backslash
\R^3$, and two nonorientable ones $M_2=\G_2\backslash \R^3,
M_3=\G_3\backslash \R^3$, also called {\em didicosm} ($c22$), {\em
first amphidicosm} ($+a2$)  and {\em second amphidicosm} ($-a2$)
respectively, in \cite{CR}. The groups $\G_i=\langle \g_{1}=B_1
L_{b_1}, \g_{2}=B_2 L_{b_2}, \Ld \rangle$, are given in the table
below, where $B_{3}=B_{1}B_{2}$, $b_{3}\equiv B_{2} b_{1} + b_{2}
\mod \Ld$, and $\Ld=\Z e_1 \oplus \Z e_2 \oplus \Z e_3$ is the
cubic lattice. All matrices $B_{i}$ are diagonal and are written
as column vectors. We indicate the translation vectors $b_i$ also
as column vectors, leaving out the coordinates that are equal to
zero.
\renewcommand{\arraystretch}{0.5}
\medskip
\begin{center}
$M_1$ \qquad
\begin{tabular}{|rc|rc|rc|}  \hline
 $B_{1}$ &  $L_{b_{1}}$  & $B_{2}$ &  $L_{b_{2}}$ & $B_{3}$ &  $L_{b_{3}}$ \\ \hline
 -1 & $\text{{\scriptsize 1/2}}$ & -1 & & 1 & $\text{{\scriptsize 1/2}}$\\
 -1 & & 1 & $\text{{\scriptsize 1/2}}$ & -1 & $\text{{\scriptsize 1/2}}$\\
  1 & $\text{{\scriptsize 1/2}}$ & -1 &  & -1 & $\text{{\scriptsize 1/2}}$\\
 \hline
\end{tabular}

\smallskip
$M_2$ \qquad
\begin{tabular}{|rc|rc|rc|}  \hline
 $B_{1}'$ &  $L_{b_{1}'}$ & $B_{2}'$ &  $L_{b_{2}'}$ & $B_{3}'$ &  $L_{b_{3}'}$ \\ \hline
 -1 & & 1 & & -1 & \\
 -1 & $\text{{\scriptsize 1/2}}$ & -1 & & 1 & $\text{{\scriptsize 1/2}}$ \\
  1 & $\text{{\scriptsize 1/2}}$ & 1 & $\text{{\scriptsize 1/2}}$ & 1 & \\
 \hline
\end{tabular}

\smallskip $M_3$ \qquad
\begin{tabular}{|rc|rc|rc|}  \hline
 $B_1'$ &  $L_{b_{1}''}$  & $B_{2}'$ &  $L_{b_{2}''}$ & $B_{3}'$ &  $L_{b_{3}''}$ \\ \hline
 -1 & & 1 & $\text{{\scriptsize 1/2}}$& -1 & $\text{{\scriptsize 1/2}}$\\
 -1 & $\text{{\scriptsize 1/2}}$& -1 & & 1 & $\text{{\scriptsize 1/2}}$\\
  1 & $\text{{\scriptsize 1/2}}$ & 1 &  & 1 & $\text{{\scriptsize 1/2}}$\\
 \hline
\end{tabular}
\end{center}

Using (\ref{dpmu}) we have
\begin{equation} \label{dpmus} \begin{split}
d_{p,\mu}(\G_1) & =  \tfrac 14\left(\tbinom{3}{p} |\Ld_\mu| + K_p^3(2)\,
(e_{\mu,\g_1}+
e_{\mu,\g_2}+e_{\mu,\g_3})\right) \\
d_{p,\mu}(\G_i) & = \tfrac 14 \left(\tbinom{3}{p} |\Ld_\mu| + K_p^3(1)\, e_{\mu,\g_1}
+ K_p^3(2)\, (e_{\mu,\g_2}+e_{\mu,\g_3})\right)
\end{split}
\end{equation}
for  $i=2,3$.

Now, to show that $M_1,M_2$ and $M_3$ are not pairwise
$p$-isospectral for any $0\le p\le 3$ we shall again use two
eigenvalues, namely those corresponding to $\mu=1$ and
$\mu=\sqrt5$.

Take $\mu=1$. Then $\Ld_1=\{\pm e_1,\pm e_2,\pm e_3\}$ and hence
$|\Ld_1|=6$. Also, $\Ld_1^{B_i}=\{\pm e_i\}$ for $1\le i\le 3$,
$\Ld_1^{B_1'}=\{\pm e_3\}$, $\Ld_1^{B_2'}=\{\pm e_1,\pm e_3\}$,
and $\Ld_1^{B_3'}=\{\pm e_2,\pm e_3\}$. For $\mu=\sqrt 5$ we see
that $\Ld_{\sqrt5}=\{ \pm(2e_i \pm e_j) : 1\le i,j\le 3 \}$, so
$|\Ld_{\sqrt5}|=24$.
 Now, we have $\Ld_{\sqrt5}^{B_i}=\Ld_{\sqrt 5}^{B_1'}=\emptyset$ for $1\le i\le 3$,
  $\Ld_{\sqrt5}^{B_2'}=\{\pm (2e_1 \pm e_3),\pm (e_1 \pm 2e_3)\}$, and
  $\Ld_{\sqrt5}^{B_3'}=\{\pm (2e_2 \pm e_3),\pm (e_2 \pm 2 e_3)\}$.

With this information one computes the following values of $e_{\mu,\g}(\G_i)$
for $1\le i \le 3$:

\smallskip
\renewcommand{\arraystretch}{0.3}
\begin{center}
\begin{tabular}{|c|ccc|ccc|} \hline &&& &&& \\
 & $e_{1,\g_1}$ & $e_{1,\g_2}$ & $e_{1,\g_3}$  & $e_{\sqrt5,\g_1}$ & $e_{\sqrt5,\g_2}$ & $e_{\sqrt5,\g_3}$ \\
& & &  & & &\\  \hline & & &  & & &\\
$M_1$ & -2 & -2 & -2  & 0 & 0 & 0 \\
$M_2$ & -2 & 0  & 0   & 0 & 0 & 0 \\
$M_3$ & -2 & 0 & -4  & 0 & 0 & -8 \\ & & &  & & & \\
\hline
\end{tabular}
\end{center}

\smallskip
\noindent
 By substituting these values back in (\ref{dpmus}) we
obtain:

\smallskip
\begin{center}
\begin{tabular}{|c||c|c|} \hline && \\
     & $d_{p,1}$  & $d_{p,\sqrt5}$ \\ && \\ \hline && \\
$\G_1$ & $\tfrac32 \big(\tbinom 3p - K^3_p(2)\big)$  &
$6\tbinom 3p$ \\ && \\
$\G_2$ & $\tfrac12 \big(3\tbinom 3p - K^3_p(2)\big)$ & $6\tbinom 3p$\\ && \\
$\G_3$ & $\tfrac12 \big(3\tbinom 3p - K^3_p(2) - 2K^3_p(1)  \big)$
 &  $6\tbinom 3p - 2K^3_p(1)$
\\ && \\ \hline
\end{tabular}
\end{center}

\smallskip
With this information and using (\ref{krawtables}), we are now in
a position to give the multiplicities for the two eigenvalues we
are considering.

\smallskip
\begin{center}
\begin{tabular}{|c|cccc|c|}  \hline & & & & &  \\
 $\mu=1$    & $d_0$ & $d_1$ & $d_2$ & $d_3$ &  $d_{f}$ \\
& & & & &\\  \hline & & & & & \\
$M_1$ & 0 & 6 & 6 & 0 & 12 \\
$M_2$ & 1 & 5 & 5 & 1 & 12 \\
$M_3$ & 0 & 4 & 6 & 2 & 12 \\ &&&&& \\ \hline
\end{tabular}
\begin{tabular}{|c|cccc|c|}  \hline & & & & &  \\
 $\mu=\sqrt 5$    & $d_0$ & $d_1$ & $d_2$ & $d_3$ &  $d_{f}$ \\
& & & & &\\  \hline & & & & & \\
$M_1$ & 6 & 18& 18& 6 & 48 \\
$M_2$ & 6 & 18& 18& 6 & 48 \\
$M_3$ & 4 & 16& 20& 8 & 48 \\ &&&&& \\ \hline
\end{tabular}
\end{center}
\smallskip
From these tables it is clear that the manifolds $M_1, M_2, M_3$
are pairwise not $p$-isospectral for any individual value of $p$,
$0\le p\le 3$.
\end{ejem}

\smallskip
\begin{ejem}
In \cite{RT} a family $\mathcal{B}_n$ of pairwise nonhomeomorphic
$\Z_2^2$-manifolds of dimension $n$  with $\beta_1=0$, is given.
This family, which
 includes Cobb's family as a rather small subfamily, grows polynomially as
  $\frac{n^5}{2^7 3^25}$. By Theorem \ref{main}, all manifolds in
  $\mathcal{B}_n$ are $\Lf$-isospectral.
\end{ejem}

\smallskip
\begin{ejem}
We first recall some facts from \cite{MRc}. Let $n$ be odd. A {\it
Hantzsche-Wendt} (or {\em HW}) {\em group} is an $n$-dimensional
orientable Bieberbach group $\G$ with holonomy group
$F\,\simeq\,\Z_2^{n-1}$ such that the action of every $B\in F$
diagonalizes on the canonical $\Z$-basis $e_1,\dots, e_n$ of
$\Ld$. The holonomy group $F$ can thus be identified to the
diagonal subgroup $\{B:Be_i = \pm e_i,
\,\,1\,\le\,i\,\le\,n,\,\,\det B = 1\}$ and $M_\Gamma
=\Gamma\backslash\R^n$ is called a {\it Hantzsche-Wendt manifold}.

Denote by $B_i$ the diagonal matrix  fixing $e_i$ and such that
$B_i e_j = -e_j$ (if $j \ne i)$, for each $1\,\le\,i\,\le\,n$.
Clearly, $F$ is generated by $B_1, B_2,\dots, B_{n-1}$.

Any HW group has the form $\G = \langle B_1
L_{b_1},\dots,B_{n-1}L_{b_{n-1}},L_\ld\,:\,\ld\,\in\, \Ld\rangle$,
for some $b_i\in\R^n$, $1\le i \le n-1$, where it may be assumed
that the components $b_{ij}$ of $b_i$ satisfy
$b_{ij}\in\{0,\frac12\}$, for $1\,\le\,i,j\,\le\,n$. Also, it is
easy to see that  ${\Ld^p(\R^n)}^{F}=0$ for any $1\le p \le n-1$,
hence all Betti numbers are 0 for $1\le p \le n-1$, thus HW
manifolds are rational homology spheres. We further recall that it
is shown in \cite{MRc} (by considering a rather small subfamily)
that the cardinality $h_n$ of the family of all HW groups under
isomorphism satisfies $h_n
>\frac {2^{n-3}}{n-1}$. Moreover, the cardinality of the pairs of
isospectral, nonisomorphic  HW groups grows exponentially
with~$n$.

All HW manifolds form a family of pairwise nonhomeomorphic compact
flat $n$-manifolds, of cardinality growing exponentially with $n$,
which by Theorem \ref{main} are mutually  $\Lf$-isospectral.
\end{ejem}

\smallskip
\begin{ejem}
Here we show that Theorem \ref{main} fails to hold without the
assumption that $F \simeq \Z_2^k$, even when the manifolds are
isospectral on functions.

(i) First, we consider the pair $M, M'$ of manifolds of dimension
6, having holonomy group $\Z_4 \times \Z_2$, studied in
\cite{MRp}, Example 5.1. Take $\tilde J:=\left[
\begin{smallmatrix} 0&1\\-1&0 \end{smallmatrix} \right]$. Let
$\G=\langle B_1L_{b_1}, B_2 L_{b_2}, \Ld \rangle$ and $\G'=\langle
B_1' L_{b'_1}, B_2' L_{b'_2}, \Ld \rangle$ where $\Ld$ is the
canonical lattice in $\R^6$ and
$$B_1=\left[ \begin{smallmatrix} \tilde J & & & \\ & \tilde J & & \\ & & 1 & \\
& & & 1 \end{smallmatrix} \right], \quad
B_2=\left[ \begin{smallmatrix} -I & & & \\ & I & & \\ & & 1 & \\
& & & 1 \end{smallmatrix} \right], \quad b_1=\tfrac{e_5}{4}, \quad
b_2=\tfrac{e_6}{2}, $$
$$B_1'=\left[ \begin{smallmatrix} \tilde J & & & &\\ & 1 & & &\\ & & -1 & &\\
& & & -1& \\ &&&&1 \end{smallmatrix} \right], \quad
B_2'=\left[ \begin{smallmatrix} -I& & & & \\ & -1& && \\ & & 1 && \\
& & & -1& \\ &&&&1 \end{smallmatrix} \right], \quad
 b_1'=\tfrac{e_6}{4}, \quad b_2'=\tfrac{e_4+e_5}{2}.$$
In \cite{MRp} it is shown that $M=\G\backslash \R^6,
M'=\G'\backslash \R^6$ are isospectral on functions --and hence
$6$-isospectral, by orientability-- but  they are not
$p$-isospectral for any $1\le p\le 5$. We claim they are not
$\Lf$-isospectral. To see this, it will be sufficient to look at
$\mu =0$. Indeed, since $d_{p,0}(\G)=\beta_p(M)$ we have that
$d_{f,0}(M)=\sum_{p=0}^6 \beta_p(M)$. The Betti numbers for $M,M'$
are given, for $0\le p \le 6$, respectively by $1,2,3,4,3,2,1$ and
$1,1,1,2,1,1,1$  (see \cite{MRd}, Remark 5.2). In this manner we
have that $d_{f,0}(M)=16$ while $d_{f,0}(M')=8$. Furthermore,  we
observe  that $d_{e,0}(M)=d_{o,0}(M)=8$ while
$d_{e,0}(M')=d_{o,0}(M')=4$. Hence $M,M'$ are not isospectral on
even forms nor on odd forms.

(ii) As a simpler example we look at a variation of (i). Indeed,
let $\G=\langle B_1L_{b_1},  \Ld \rangle$ and $\G'=\langle B_1'
L_{b'_1}, \Ld \rangle$ with $B_j, b_j$ as in (i), $j=1,2$. Then
$\Gamma, \Gamma'$ have holonomy group $\Z_4$ and a computation
shows that the Betti numbers are, in this case, respectively given
by $1,2,5,8,5,2,1$ and $1,2,3,4,3,2,1$. For instance, to verify
the values in the case of $\beta_2$, we note that a basis of
vectors fixed by $B_1$ on $\Lambda^2(\R^6)$ is given by $e_1\wedge
e_2,  e_3\wedge e_4,e_5\wedge e_6, e_1\wedge e_4,e_2\wedge e_3$
while for $B'_1$ a basis of  fixed vectors is  $e_1\wedge e_2,
e_3\wedge e_6,e_4\wedge e_5$.

Thus, in this case we  have $d_{f,0}(M)=24$ while
$d_{f,0}(M')=16$, hence $M$ and $M'$ are not isospectral on forms.
\end{ejem}

\section{Large families of manifolds isospectral on forms}
In this section we will exhibit large families of
$\Z_2^k$-manifolds, pairwise nonhomeomorphic to each other, which
by Theorem \ref{main} will be isospectral on forms. In particular,
for each $n$, we shall construct a family of $n$-dimensional
$\Z_2^{n-1}$-manifolds with cardinality of order approximately
\begin{small}$(\sqrt{2})$\end{small}$^{^{n^2}}$, for $n$ large.
Consider the subgroup of $I(\R^n)$ with set of generators $\{
C_{i}L_{c_{i}}, L_{e_j} : 1 \leq i \leq n-1, \, 1 \le j \le n \}$,
where $C_i:=\text{diag}(1,\dots,1,\underset i{-1},1,\dots,1)$, and
$c_{i} := \frac{e_{i+1}}2 + \sum_{j=1}^i c_{ji} e_j$ for some
choice of $c_{ji} \in \{0, \frac 12\}$. Here $\{e_{1},e_{2}, \dots
,e_{n}\}$ denotes the canonical basis of $\Ld=\Z^{n}$. Similarly
as in Example \ref{ej3.2} above, we show in (\ref{LeeScz}) such a
group in {\it column notation}, placing the coordinates of the
translation vectors $c_i$ as subindices in each column.

\renewcommand{\arraystretch}{0.35}
\begin{equation}\label{LeeScz}
\begin{array}{|cccccc|}
\hline
{\phantom{_{8_{8_8}}}}C_1^{\phantom{{8^8}^8}} & C_2 & C_3 & \dots & C_{n-1} &  \bar{C}_n \\
\hline
-1^{\phantom 8}_{\phantom{\frac 12}} & {\phantom -}1_* & {\phantom -}1_* & \dots & {\phantom -}1_* & -1_* \\
{\phantom -}1_{\frac 12} & -1{\phantom 8} & {\phantom -}1_* & \cdots& {\phantom -}1_* & -1_* \\
1 & {\phantom -}1_{\frac 12} & -1{\phantom 8} & \ddots & \vdots & \vdots \\
1 & 1 & {\phantom -}1_{\frac 12} & \ddots & {\phantom -}1_* & -1_* \\
\vdots & \vdots & \vdots & \ddots & -1_{\phantom{\frac 12}} & -1_* \\
1 & 1 & 1 & \dots & {\phantom -}1_{\frac 12} & {\phantom
-}1_{\frac 12} \\[.2cm]
\hline
\end{array}
\end{equation}
where each $*$ is $0$ or $\frac 12$, depending on the choice of
the $c_{ji}$'s. We have added an extra column $ \bar C_n$
corresponding to the product $\bar C_n:=C_1C_2\dots C_{n-1}$ and
we take the respective column vector  $c_n\equiv
c_1+\dots+c_{n-1}$ mod~$\Ld$ and  having coordinates in $\{0,\frac
12\}$.

\begin{figure}[!htb]
$\begin{array}{|cccccc|} \hline
0 & 0 & \cdots & \cdots & 0 & 0 \\
\frac 12 & 0 & \cdots & \cdots & 0 & \frac12 \\
0 & \frac 12 & \ddots & \ddots& \vdots & \vdots \\
0 & 0 & \ddots & \ddots & \vdots & \vdots \\
\vdots & \vdots & \ddots & \ddots & 0 & \frac 12 \\
0 & 0 & \cdots & 0 & \frac 12 & \frac 12 \\
\hline
\end{array}$
\quad\quad \raisebox{-.6cm}[0pt][0pt]{\input{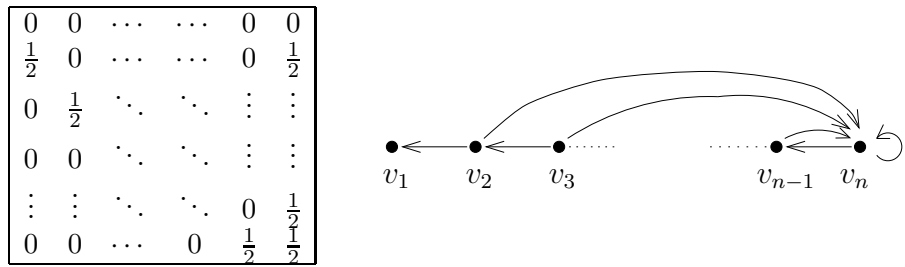}}
\caption{The matrix and the graph of the group $K_n$.}
\label{secgraph}
\end{figure}

In (\ref{LeeScz}), in the case when all $*$'s in the first $n-1$
columns equal zero (and thus the $*$'s in the $n$-th column are
$\frac 12$'s, except for the one  in the entry $(1,n)$ which is
zero), the corresponding group, which we will denote by $K_n$ (see
Figure~\ref{secgraph}), was introduced in \cite{LS} and is known
to be torsion-free, i.e., a Bieberbach group. Here, we will prove
that this is true in the more general case above. We shall denote
by $\mathcal{K}_n$ the family consisting of all groups constructed
in this manner.

\begin{prop} All groups in~$\mathcal{K}_n$ are Bieberbach groups.
\end{prop}
\begin{proof}
It is clear that the groups are Euclidean crystallographic groups.
Hence, we need only  show that they are torsion-free. Every
element in each group is either a translation $L_\ld$ with
$\ld\in\Ld$, or it is of the form
\begin{equation}\label{prodgen}
\g:=C_{i_1}L_{c_{i_1}} \dots C_{i_k}L_{c_{i_k}} L_\ld, \quad
i_1<\dots<i_k,\quad 1\le k\le n-1.
\end{equation}

The translations $L_\ld$, $\ld\not =0$, are clearly not elements
of finite order. Concerning the remaining elements, we observe
that on the $(i_k+1)$-th coordinate the product $\g$ in
(\ref{prodgen}) acts as the translation by $\frac 12+\ld_{i_k+1}$,
with $\ld_{i_k+1}\in\Z$ and $\ld=(\ld_1,\dots,\ld_n)$. This is so
since the action of $C_{i_j}$ on $e_{i_k+1}$ is trivial for
$j=1,2,\dots,k$, and the $(i_k+1)$-th coordinates of the
translational parts $c_{i_j}$ are equal to zero, for
$j=1,\dots,k-1$, and $\frac 12$ for $j=k$.

Now, $\g^2$ is a translation and its $(i_k+1)$-th coordinate is
$2(\frac 12+\ld_{i_k+1})=2\ld_{i_k+1}+1$. In general, the
$(i_k+1)$-th coordinate of the $m$-th power of $\g$ in
(\ref{prodgen}) equals $m(\frac 12+\ld_{i_k+1})$, hence $\g^m\ne
Id$ for every $m\ne0$.
\end{proof}

\begin{rem}
Manifolds of dimension $n$ and holonomy group $\Z_2^{n-1}$ have
been called in \cite{RS} {\it generalized Hantzsche-Wendt
manifolds}, or GHW manifolds for short. They necessarily have
diagonal holonomy representation. The family ${\mathcal K}_n$ is
properly contained in this larger class. It is not difficult to
see that the holonomy representation in (\ref{LeeScz}) is the only
possible one for GHW manifolds with first Betti number one (see
\cite{RS}).
\end{rem}

\

Next we shall show that all the Bieberbach groups in $\mathcal
K_n$ are pairwise not isomorphic. In \cite{MRc}, we have attached
a directed graph to any orientable GHW manifold, with
diffeomorphism of manifolds corresponding to isomorphism of
graphs, essentially. It is also possible  to do the same for
arbitrary GHW manifolds, i.e., to associate a directed graph with
$n$ vertices to any $n$-dimensional GHW manifold. We will do this
in the case of the family ${\mathcal K}_n$. This graph will be
helpful to better understand the elements in our family.

Firstly, we replace the array in (\ref{LeeScz}) by an $n\times n$
array $A$ of $0$'s and $\frac 12$'s by keeping just the
translational parts mod $\Z^n$. We observe that the total number
of $\frac 12$'s in each row must be even, since the last column in
$A$ is the sum mod $\Z^n$ of the others. Thus

\renewcommand{\arraystretch}{1}
\begin{equation}\label{array}
A=
\begin{array}{|cccccc|}
\hline
0 & * & * & \cdots & * & * \\
\frac 12 & 0 & * & \cdots & * & * \\
0 & \frac 12 & 0 & \ddots & \vdots & \vdots \\
0 & 0 & \frac 12 & \ddots & * & \vdots \\
\vdots & \vdots & \vdots & \ddots & 0 & * \\
0 & 0 & 0 & \cdots & \frac 12 & \frac 12 \\
\hline
\end{array}
\end{equation}
where each $*$ can be equal to $0$ or $\frac 12$.

We recall that these arrays are in a one-to-one correspondence
with groups in $\mathcal K_n$. We will associate to each such
array (or group in $\mathcal{K}_n$) a directed graph having a
fixed set of $n$ vertices $\{v_1,\dots,v_n\}$ and so that there is
an arrow issuing from vertex $v_i$ into vertex $v_j$ if and only
if the entry $(i,j)$ in $A$ equals~$\frac 12$.

\begin{rem}
Note that similar matrices $A$ were used in \cite{MRc} to describe
orientable GHW manifolds (there called  HW manifolds), however the
translational parts, shown as columns of the matrix $A$ here, are
shown as rows in \cite{MRc}. We also observe that another option
to define these graphs would have been to  `colour' vertex $v_n$
leaving out the arrow joining $v_n$ to itself. In \cite{MRc} these
loops were omitted since all the vertices were of the same kind.
\end{rem}

To illustrate the definition above, we display some GHW groups and
their graphs.

In Figure \ref{secgraph}, we show the array and the graph
corresponding to the group $K_n$. Note that the arrows going from
right to left in the figure of the graph will be present in every
graph corresponding to a group in $\mathcal K_n$.

Next, we will show the graphs of GHW groups in dimensions 2 and 3.
In dimension 2, the Klein bottle group (which is isomorphic to
$K_2$) belongs to~$\mathcal K_2$.

\begin{figure}[!htb]
$\begin{array}{|cc|} \hline
0 & 0  \\
\frac 12 & \frac 12 \\
\hline
\end{array}$
\qquad\raisebox{-.3cm}[0pt][0pt]{\mbox{\includegraphics*{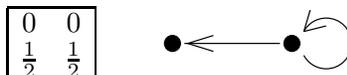}}}
\caption{The graph of the Klein bottle.} \label{kleinfig}
\end{figure}

\begin{figure}[!htb]
\begin{tabular}{cc}
$\begin{array}{|ccc|} \hline
0 & 0 & 0 \\
\frac 12 & 0 & \frac 12 \\
0 & \frac 12 &  \frac 12 \\
\hline
\end{array}$ & \qquad
\raisebox{-.3cm}[0pt][0pt]{\mbox{\includegraphics*{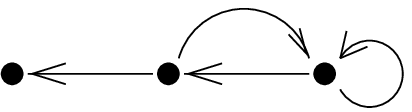}}}
\\
 & \\
  & \\
$\begin{array}{|ccc|} \hline
0 & \frac 12 & \frac 12 \\
\frac 12 & 0 & \frac 12 \\
0 & \frac 12 &  \frac 12 \\
\hline
\end{array}$ & \qquad
\raisebox{-.7cm}[0pt][0pt]{\mbox{\includegraphics*{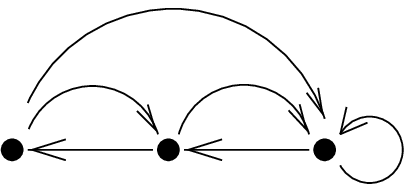}}}
\end{tabular}
\caption{The graphs of the first and second amphidicosms $+a2$ and
$-a2$.} \label{a2pmfig}
\end{figure}

\begin{figure}[!htb]
\centerline{\mbox{\includegraphics*[scale=.7]{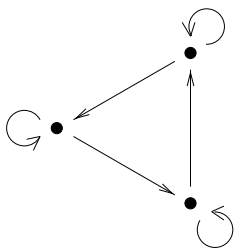}}}
\caption{The graph of the didicosm $c22$.} \label{c22pfig}
\end{figure}

In dimension 3, out the three existing GHW groups, two of them,
the first amphidicosm $+a2$ and the second amphidicosm $-a2$,
belong to $\mathcal K_3$ while the other one, the didicosm $c22$
(or Hantzsche-Wendt manifold), does not (see figures \ref{a2pmfig}
and \ref{c22pfig}).

In dimension 4, there are twelve GHW manifolds (see \cite{CS} for
instance), ten having first Betti number equal to one, i.e.\@
$\beta_1=1$ (see figures 5 and 6) and two having $\beta_1=0$. Out
of these, eight manifolds are in $\mathcal K_4$. They are given by
the array

\begin{equation}\label{arraydim4}
\begin{array}{|cccc|}
\hline
0 & x & y & x+y \\
\frac 12 & 0 & z & z+ \frac 12 \\
0 & \frac 12 & 0 & \frac 12 \\
0 & 0 & \frac 12 &  \frac 12 \\
\hline
\end{array}
\end{equation}
\noindent where $x,y,z \in \{0, \frac 12\}$ and the sums are taken
mod $\Z$. If we choose the eight different possibilities for
$x,y,z\in\{0,\frac 12\}$, we obtain the eight groups in $\mathcal
K_4$. The eight graphs corresponding to these manifolds have some
common features, as can be seen in Figure 5.

\begin{figure}[!htb]\label{ksfig}
\begin{tabular}{r@{\hskip2cm}ccc}
 & $x$ & $y$ & $z$
\\
{\includegraphics*[scale=.8]{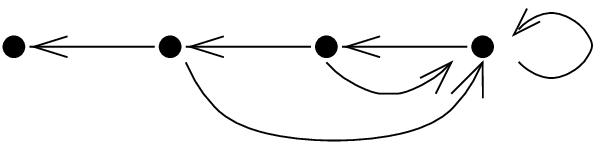}} & 0 & 0 & 0
\\
{\includegraphics*[scale=.8]{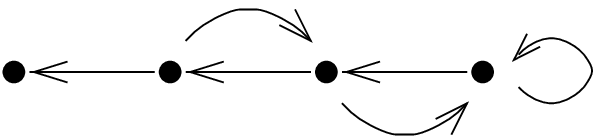}} & 0 & 0 & $\frac 12$
\\
{\includegraphics*[scale=.8]{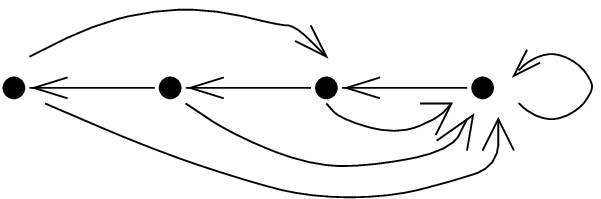}} & 0 & $\frac 12$ & 0
\\
{\includegraphics*[scale=.8]{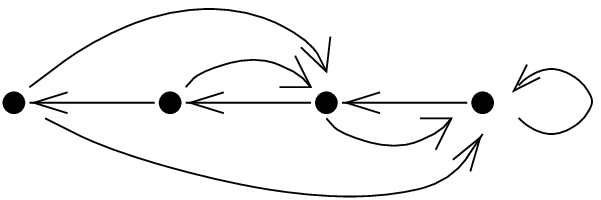}} & 0 & $\frac12$ & $\frac 12$
\\
{\includegraphics*[scale=.8]{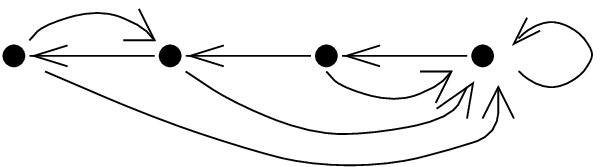}} & $\frac12$ & 0 & 0
\\
{\includegraphics*[scale=.8]{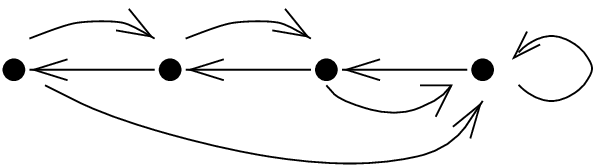}} & $\frac12$ & 0 & $\frac
12$
\\
{\includegraphics*[scale=.8]{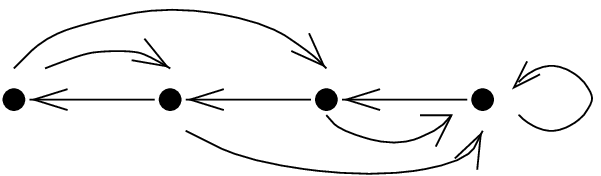}} & $\frac12$ & $\frac 12$
& 0
\\
{\includegraphics*[scale=.8]{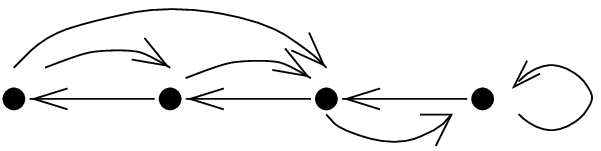}} & $\frac12$ & $\frac 12$
& $\frac 12$
\end{tabular}

\caption{The graphs of the eight manifolds in ${\mathcal K}_4$. On
the right, we show the values of $x,y,z$ of the corresponding
arrays.}
\end{figure}

\begin{figure}[!htb]\label{k910}
{\includegraphics*[scale=.9]{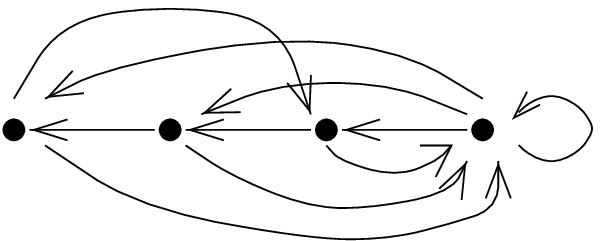}} \hskip1cm
{\includegraphics*[scale=.9]{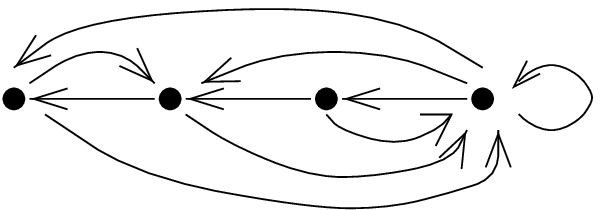}} \caption{The graphs of
the two 4-dimensional GHW manifolds with $\beta_1=1$ which are not
in ${\mathcal K}_4$.}
\end{figure}

\begin{prop} Let $\Gamma$ and $\Gamma'$ be groups in ${\mathcal K}_n$
corresponding to arrays $A$ and $A'$ as in (\ref{array}),
respectively. Then $\Gamma\,\simeq\,\Gamma'$ if and only if
$A=A'$, or equivalently, if and only if their associated graphs
$G$ and $G'$ are isomorphic.
\end{prop}

\begin{proof}
We first check that if two directed graphs, $G,G'$, attached to
arrays as in (\ref{array}) are isomorphic (as directed graphs)
then $G=G'$. We will do this by showing that each vertex $v_i$,
$1\le i\le n$, is completely determined by the isomorphism class
of the graph. The vertex $v_n$ is determined since it is the only
vertex with an arrow going to itself. Now, there is only one more
arrow issuing from $v_n$ and it goes to $v_{n-1}$, so this
determines $v_{n-1}$. Also, there is only one arrow issuing from
$v_{n-1}$ and going to vertices different from $v_n$. This arrow
goes to $v_{n-2}$, so this determines $v_{n-2}$. Continuing in
this  way we get  $v_2$  determined; furthermore $v_2$ is the only
vertex which has an arrow  going to $v_1$, hence determining
$v_1$. Thus, all the vertices are determined. In other words,  an
isomorphism $\phi$ between two of these graphs $\phi:G\rightarrow
G'$ (recall that the set of vertices is the same in both graphs)
satisfies $\phi(v_i)=v_i$ for every $i=1,2,\dots,n$, hence
$\phi=Id$ and thus $G=G'$.

By the definition of the graph, one has that $G=G'$ if and only $A=A'$.
Thus, we will be done if we prove that $\G\simeq\G'$ implies $A=A'$.

Let $\G' = \langle C_i L_{c'_i} : \,1\,\le\,i\,\le\,n-1;
L_\Lambda\rangle$. By Bieberbach's second theorem, an isomorphism
between $\Gamma$ and $\Gamma'$ must be given by conjugation by an
affine motion $\delta = D L_d$, $D\,\in\,\text{GL}_n(\R)$,
$d\in\R^n$, i.e. $\Gamma' = \delta \Gamma \delta^{-1}$. This
implies, in particular, that
\begin{equation}\label{conjugation}
\delta C_i L_{c_i} \delta^{-1} = D C_i D^{-1} L_{D(c_i+(C_i-Id)
d)}.
\end{equation}

Since $\bar C_n$ is the only matrix in the holonomy group with
exactly $n-1$ diagonal elements equal to $-1$, we must have $D\bar
C_nD^{-1}=\bar C_n$. Similarly, there is a permutation $\sigma\in
S_{n-1}$ such that $DC_iD^{-1}=C_{\sigma(i)}$ for $1\le i\le n-1$.
Also,
\begin{equation}\label{modL}
c'_{\sigma(i)} \equiv D(c_i + (C_i-Id)d)\mod\Lambda.
\end{equation}

Now we take into account that $c_i$ and $c_i'$ are of the form
$\sum_{j=1}^{i+1} \epsilon_{ji} e_j$ where $\epsilon_{ii}=0$ and
$\epsilon_{i+1\,i}=\frac 12$ for $i=1,2,\dots,n-1$,
$\epsilon_{ji}=0$ or $\frac 12$ for $j<i$ and $\epsilon_{nn}=\frac
12$. We note that equation (\ref{modL}) allows to change modulo
$\Z$ only the coordinates in $c_i$ in which $C_i$ acts by $-1$,
while the other coordinates cannot change (modulo $\Z$). In
particular, since the entries $(ji)$ with $j>i$ in $A$ correspond
to the action of $C_i$ as the identity, they remain unchanged
after conjugation by $\delta$. Then we see that $C_{n-1}$ must be
matched (via the isomorphism) to $C'_{n-1}$ (since they are the
only $C_i$'s with $\frac 12$ in the $n^{\rm th}$ coordinates of
their translation vectors). This implies that $\sigma(n-1)=n-1$.
In the same way, we see that for each $i$, $1\le i\le n-2$, $C_i$
must be matched with $C'_i$, thus $\sigma(i)=i$ for $1\le i\le
n-1$. Therefore, it follows that the permutation $\sigma$ must be
the identity.

This implies that $D$ is diagonal with eigenvalues $\pm 1$.  By
taking into account that we have chosen the coefficients in the
main diagonal in $A$ and in $A'$ to be zero for the first $n-1$
entries, we must have $d\equiv 0\mod \Z$. Hence, conjugation by
$\delta = D L_d$ produces an automorphism of $\G$. Thus, $\G=\G'$,
and hence we have $A=A'$ for the corresponding arrays, which
completes the proof.
\end{proof}

Now, it is easy to compute the cardinality of $\mathcal{K}_n$,
since there are two choices for each  entry $(i,j)$ with $1\le
i<j<n$:

\begin{coro}\label{cardkn}
There are $2^{\frac{(n-1)(n-2)}2}$ Bieberbach groups in
$\mathcal{K}_n$, all of them pairwise nonisomorphic to each other.
\end{coro}

If we put this corollary together with Theorem \ref{main}, we
have:

\begin{coro}\label{isoskn}
There exists a family of $2^{\frac{(n-1)(n-2)}2}$ compact
$n$-manifolds, isospectral on forms and  pairwise
non\-homeomorphic to each other.
\end{coro}

\begin{rem}
(i) It is easy to see that there are larger families with similar
properties as $\mathcal{K}_n$. Indeed, in \cite{RS} it was shown
that, for a given $n$, there are $[(n+1)/2]$ different integral
holonomy representations for GHW manifolds. For each of these
representations, one can define  a family of Bieberbach groups in
a similar way as for $\mathcal{K}_n$ above, and all the resulting
flat manifolds will be isospectral on forms by Theorem \ref{main},
yet pairwise nonhomeomorphic. Thus, this procedure should allow to
multiply the number in Corollaries \ref{cardkn} and \ref{isoskn}
by a factor of $[(n+1)/2]$, approximately. However, this does not
improve the result significantly.

(ii) If one considers families of Bieberbach groups with holonomy
group $\Z_2^k$ for some  $k$ with $\frac n2\le k<n-1$, one should
obtain larger families of manifolds than in the case $k=n-1$,
pairwise nonhomeomorphic to each other and again isospectral on
forms. A support for this claim is given by the classification of
low dimensional Bieberbach groups (see \cite{CS}).

(iii) The manifolds in ${\mathcal K}_n$ are all nonorientable. By
using a {\em duplication method} (see for instance \cite{BDM})
applied to the manifolds in ${\mathcal K}_n$ one obtains
$2^{\frac{(n-1)(n-2)}2}$, orientable, nonhomeomorphic manifolds of
dimension $2n$ isospectral on forms.

(iv) Using the methods in \cite{MRd} (see Thm.\ 3.12) and in
\cite{MRl} (see Prop.\ 4.7) one can show that the manifolds in
$\mathcal{K}_n$ are, generically, not $p$-isospectral for any
value of $p$, $0\le p\le n$.

\end{rem}



\end{document}